\documentclass{beamer}
\usepackage[T1]{fontenc} 
\usepackage[normalem]{ulem}
\usepackage{amsmath}
\usepackage{amsfonts}
\usepackage{amsthm}
\usepackage{verbatim}
\usepackage{geometry}
\usepackage{graphicx}
\usepackage{movie15}
\usepackage{animate}
\usepackage{ulem}
\usepackage{color}

\newsavebox\dummy

\definecolor{VF}{rgb}{0.15,0.7,0.15}

\title{Optimal (static output feedback) design through Direct Search methods}
\author{Emile Simon\\UCL Belgium}
\date{December 12, 2011\\ $50^{th}$ CDC-ECC}
\begin{document}

\begin{frame}
\maketitle
\end{frame}

\begin{frame}{Aim of this work:}

Putting forward a class of optimization methods\\ 
largely overlooked in systems and control:\\ 
the direct search methods.\\

\vspace{1cm}

They can however be very adequate and powerful\\
for many complex problems (possibly including yours).  
  
\setcounter{framenumber}{1}
  
\end{frame}


\begin{frame}
\section{Direct search optimization}
\tableofcontents[sections={1-3}, currentsection, hideothersubsections]
\end{frame}

\begin{frame}{Direct search optimization methods}
\begin{itemize}
  \item What are these methods?
  \vspace{0.4cm}
  \item Why are they overlooked in systems and control?
  \vspace{0.4cm}
  \item Why should they be used?
  \end{itemize}
  
 \setcounter{framenumber}{2}
  
\end{frame}

\begin{frame}{What are these methods?}

\begin{itemize}
 \item The aim is to min(/max)imize an objective function $f(x):\mathbb{R}^n\rightarrow \mathbb{R}$, possibly under constraints,\\
 starting from one (or several) feasible initial solution(s).
  \vspace{0.4cm}
  \item \textit{Only use function evaluations} to decide how to explore $\mathbb{R}^n$,
 \vspace{0.4cm}
 \item \textit{need no gradient or Hessian information}.
 \vspace{0.4cm}
 \item In this talk, direct search = derivative-free optimization
  \vspace{0.4cm}
 \item See e.g. webpages of Luis Vicente and Charles Audet.
 \end{itemize}
    
\end{frame}


\begin{frame}{Why are they overlooked in systems and control?: A history}
\begin{itemize}
  \item 1957-1961: The birth (Box, Davidon, Hooke, Jeeves) 
  \vspace{0.2cm}
  \item 1961-1971: The golden age (1965 Nelder-Mead)
  \begin{itemize}
  \item Very efficient in practice.
  \end{itemize}
  \vspace{0.2cm}
  \item 1972-1989: The downfall (after survey of W. Swann 1972)
  \begin{itemize}
  \item No proofs of convergence and can sometimes be slow.
  \end{itemize}
  \vspace{0.2cm}
  \item 1990-...: The resurrection = Proofs of convergence
  \begin{itemize}
  \item Torczon and Dennis: MDS, proof on smooth $f(x)$: 1997
  \item Audet and Dennis: MADS, proof on non-smooth $f(x)$: 2003
  \item Vicente and Custodio: proof on discontinuous $f(x)$: 2010
  \end{itemize}
  \end{itemize}
  
  Last item \textit{arrived too late for systems and control}:\\
  From the early nineties, our field got saturated by convex optimization and LMIs, solvable very efficiently but...
  
\end{frame}

\begin{frame}{Why should they be used? (1/2)}
\begin{itemize}
 \item Current problems of interest have \textit{non-convex} feasible sets.
\vspace{0.4cm}
 \item The best that can be done is to propose methods with\\ 
\textit{guaranteed convergence to locally optimal solutions.}
\vspace{0.4cm}
 \item There exist (new) direct search methods \textit{guaranteed to converge even on most non-smooth or discontinuous $f(x)$}.
  \end{itemize}
  
\end{frame}

\begin{frame}{Why should they be used? (2/2)}

When (clean) gradients of $f(x)$ are not (easily) available.\\
\vspace{0.4cm}
The main advantages are:
\vspace{0.1cm}
\begin{itemize}
	\item Recent strong convergence guarantees
	\vspace{0.3cm}
	\item Exponential performance of computers
	\vspace{0.3cm}
	\item Simulations are more routine and accurate
	\vspace{0.3cm}
	\item Ease of use and implementation
	
\end{itemize}
   
\end{frame}

\begin{frame}{Optimization methods: A comparison}

A basic comparison of optimization methods:

\begin{table}[h]
\begin{flushleft}

\begin{tabular}{|l||c|c|}
\hline
Class of methods  & Computational time & Adequate problems \\
\hline
\hline
Convex/LMI & Very efficient\footnote{May need to introduce a lot of additional variables to convexify\\ \vspace{0.4cm}} & Many specific \\
 & & (sub)problems \\
\hline
Gradient-based & Efficient & Clean derivatives  \\
 & & must be available\\
\hline
Derivative-free & May be slow  &`Any' $f(x):\mathbb{R}^n\rightarrow \mathbb{R}$\\
 &  & ($n<25-100$)\\
\hline
\end{tabular}
\end{flushleft}
\end{table}

\end{frame}

\begin{frame}
\section{Problems in the paper (SOF)}
\tableofcontents[sections={1-3}, currentsection, hideothersubsections]
\addtocounter{framenumber}{-1}
\end{frame}

\begin{frame}{The problems considered in the paper (1/2)}

Static Output Feedback optimization, $m\times p$ MIMO LTI systems:\\
\vspace{0.3cm}
1) Stabilization of the closed-loop:\\
\begin{center}
Find $K \in \mathbb{R}^{m\times p}$ s.t. $\sigma(A+BKC)<0$  (continuous time)\\
\end{center}
2 and 3) Min. the $\mathcal{H}_2$ or $\mathcal{H}_\infty$ norm of the performance channel:\\
\begin{center}
$\min_K ||T_{wz}(K,s)||_2$ or $\min_K ||T_{wz}(K,s)||_\infty$, s.t. $\sigma(A+BKC)<0$\\
\end{center}
\vspace{0.2cm}
With unstable models from COMPlib : a library of actual and academic models, currently often used for benchmarking.

\end{frame}

\begin{frame}{The problems considered in the paper (2/2)}

Non-convex: \sout{LMIs} (and iterative LMI schemes rarely converge)\\
\vspace{0.4cm}
The objectives admit however \textit{gradients} or \textit{Clarke subgradients},\\
which should then be used for optimization purposes.\\
\vspace{0.4cm}
This is implemented in the methods HIFOO and \texttt{hinfstruct}.\\
\vspace{0.4cm}
Paper's benchmarks: comparison of the objectives values reached and computational times required, by a DS method and by HIFOO.\\
\vspace{0.4cm}
Motivation: to verify convergence and assess performance of DS.\\
Then other $f(x)$ \textit{without gradients} can be considered.

\end{frame}

\begin{frame}{Note on the direct search method used}

Currently, the best two derivative-free methods are certainly MADS and SID-PSM, having strong theoretical convergence guarantees.\\
\vspace{0.4cm}
In practice, a very efficient method easily accessible and usable is the Nelder-Mead algorithm \textit{restarted at the last solution until no improvement is obtained} (to a given accuracy).\\
\vspace{0.4cm}
This generates a set of dense (enough) exploring directions in the search space $\mathbb{R}^n$ (but without formal theoretical guarantee).\\
\vspace{0.4cm}
Packs a serious punch in practice, sufficing to make the point.

\end{frame}

\begin{frame}{Summary of 13000 tests}

The direct search method converges as expected, and is reasonably fast (meant for off-line design), despite not using gradients:
\vspace{0.2cm}
\begin{itemize}
\item 1) \textit{Stabilization success:} DS 92.3\%, HIFOO 90.6\%.\\ 
Ratio computational time DS on HIFOO $\cong$ 0.74
\vspace{0.3cm}
\item 2) \textit{Performance channel $\mathcal{H}_2$ norm minimization:}\\ 
Similar: 60\%	/ Better for DS: 14\%	/	Better for HIFOO: 26\% \\
Ratio computational time DS on HIFOO $\cong$ 10
\vspace{0.3cm}
\item 3) \textit{Performance channel $\mathcal{H}_\infty$ norm minimization:}\\ 
Similar: 40\%	/ Better for DS: 23\%	/	Better for HIFOO: 37\% \\
Ratio computational time DS on HIFOO $\cong$ 1.3
\end{itemize}

\end{frame}

\begin{frame}
\section{Other problem, without gradient}
\tableofcontents[sections={1-3}, currentsection, hideothersubsections]
\addtocounter{framenumber}{-1}
\end{frame}

\begin{frame}{Time-response shaping problem}

Objective considered: to optimize explicitly the features of a time response $z(x,t)$, with the controller parameters $x\in \mathbb{R}^n$.\\
\vspace{0.5cm}
Interesting objective: $f(x)$ = \textcolor{blue}{$t_r$} + $\lambda$ \textcolor{VF}{$max_{dev}$}, where:\\
\vspace{0.4cm}
\textcolor{blue}{$t_r$}$ \in \mathbb{R}^+$ is the \textcolor{blue}{\textit{rise time}} needed by the response $z(x,t)$ to reach a desired settling region (above $z_{max}(t)$ and under $z_{min}(t)$),\\
\vspace{0.4cm}
\textcolor{VF}{$max_{dev}$}$ \in \mathbb{R}^+$ is the \textit{maximum deviation} of the response $z(x,t)$ outside of the desired region (for $t>0$ above, and $t>t_r$ under)\\ 
\vspace{0.5cm}
and $\lambda \in \mathbb{R}^+$ is the scalar trade-off parameter.

\end{frame}

\begin{frame}{Difficulty and solution}

The objective \textcolor{VF}{$max_{dev}$} is non-smooth but locally Lipschitz: 
\vspace{0.1cm}
\begin{itemize}
\item Subgradients can be used for minimization to locally optimal solutions (see Bompart, Apkarian and Noll 2008).
\vspace{0.3cm}
\item (New) DS methods also converge on \textcolor{VF}{$max_{dev}$} (but slower).
\end{itemize}
\vspace{0.6cm}

The objective \textcolor{blue}{$t_r$} is however more difficult: 
\vspace{0.1cm}
\begin{itemize}
\item It is actually discontinuous and no gradient analysis exists.\\
\vspace{0.3cm}
\item (New) DS have guarantees of convergence for such objectives!\\
\end{itemize}

\end{frame}

\begin{frame}{Two results}

Time-response $z(x,t)$ shaping examples:\\
\vspace{0.2cm}
Optimization of PID parameters ($x=[K_p,K_i,K_d]$), by $\min_x f(x)$.\\
\vspace{0.2cm}
Choice of time envelope: $z_{min}(t)=0.98$ and $z_{max}(t)=1.02$\\
\vspace{0.4cm}
Two animated figures present two optimizations.

\end{frame}

\begin{frame}{Results animation (1/2)}
Initial solution: Ziegler-Nichols parameters. 

\begin{center}
\animategraphics[controls,width=3in]{10}{Film7Quat/film_}{1}{209}\\

(click on figure to launch animation)
\end{center}

\addtocounter{framenumber}{-1}

\end{frame}

\begin{frame}{Results animation (2/2)}

Initial solution: random parameters, not stabilizing.

\begin{center}
\animategraphics[controls,width=3in]{12}{Film8Quat/film_}{1}{320}\\

(click on figure to launch animation)
\end{center}

\end{frame}

\begin{frame}{Take-home messages}

\begin{itemize}
  \item There definitely exists alternatives to LMIs.\\ Which method is adequate for the considered problem?
 \vspace{0.4cm}
  \item Direct search methods now have strong theoretical convergence guarantees, and are acceptable in practice (efficient computers, methods and simulators).
  \vspace{0.4cm}
  \item Should be much more considered for many open\\ problems of optimization in systems and control.
	                
\end{itemize}

\end{frame}

\begin{frame}
\addtocounter{framenumber}{-1}

  \frametitle{} 
  
    \begin{flushleft}
			{Thank you for your attention!}
    \end{flushleft}
    
    \vspace{1.5cm}
    
    \begin{center}
      {Questions?}
    \end{center}
    
    \vspace{1.5cm}   
    
    \begin{flushright}
    	{Emile.Simon@uclouvain.be}
    \end{flushright}
    
\end{frame}

\end{document}